\documentclass[12pt]{article}
\usepackage{amssymb,latexsym,amsmath,euscript}\usepackage{amsfonts}
\usepackage[cp1251]{inputenc}\usepackage[epsf]{}\usepackage[english,russian]{babel}\usepackage{graphicx}
\newtheorem{theorem}{Theorem}
\newtheorem{lemma}{Lemma}
\newtheorem{proposition}{Proposition}[section]

\newtheorem{conjecture}{Conjecture}

\begin{document}

\title{Nonsingular structural stable chaotic 3-flows of attractor-repeller type}
\author{Zhentao Lai\and V. Medvedev\and Bin Yu\and E. Zhuzhoma}
\date{}
\maketitle

\renewcommand{\abstractname}{Absrtact}
\renewcommand{\refname}{Bibliography}
\renewcommand{\figurename}{Fig.}

\begin{abstract}
We show that any orientable closed 3-manifold $M$ admits structurally stable non-singular flow $f^t$ whose non-wandering set $NW(f^t)$ consists of a 2-dimensional expanding attractor and finitely many repelling periodic trajectories. For $M=\mathbb{S}^3$, we prove that the set of repelling periodic trajectories can be an arbitrary link provided that this link contains the figure eight knot. When a link consists of a unique repelling periodic trajectory (not necessarily a figure eight knot), we prove that this trajectory cannot be a torus knot. For any closed \(3\)-manifold \(M\),
we show that there does not admit any structurally stable non-singular flow $f^t$ whose non-wandering set $NW(f^t)$ consists of a 2-dimensional expanding attractor and a repelling periodic trajectory so that the repelling periodic trajectory is a trivial knot (i.e., it bounds a disk in \(M\)).
\end{abstract}

\section*{Introduction}

A dynamical system is considered \textit{chaotic} if its topological entropy is positive.
A prominent example of a chaotic flow is the geodesic flow on a closed hyperbolic surface with negative curvature.
Anosov \cite{Anosov67} proved that geodesic flows on closed Riemannian manifolds with negative curvature are rough in the sense of Andronov-Pontryagin \cite{AndronovPontryagin1937}.
Consequently, such geodesic flows are structurally stable within the class of $C^1$-smooth dynamical systems.
Anosov flows are a prime example of such structurally stable chaotic flows.

Another example of a structurally stable, nonsingular chaotic 3-flow is the dynamical suspension over the well-known Smale horseshoe \cite{Smale67}.
For this construction, let $\mathbb{S}^n$ denote the $n$-sphere, where $n \geq 1$.
The original Smale horseshoe construction can be extended from the unit square to the 2-sphere $\mathbb{S}^2$, yielding a structurally stable diffeomorphism $f: \mathbb{S}^2 \to \mathbb{S}^2$.
This diffeomorphism possesses a zero-dimensional, nontrivial, non-wandering set that is locally homeomorphic to the product of two Cantor sets.
The resulting dynamical suspension, denoted $sus^t(f)$, is a structurally stable nonsingular flow on $\mathbb{S}^2 \times \mathbb{S}^1$.
This flow features a one-dimensional, nontrivial, saddle-type non-wandering set.

According to \cite{Hayashi97,Mane88b,Rob76}, any structurally stable dynamical system satisfies Smale's Axiom A \cite{Smale67}.
This axiom stipulates that the non-wandering set is hyperbolic and that the set of periodic orbits is dense within it.
For fundamental concepts in Dynamical Systems, see \cite{grinmedpoch-book-2016-eng,GrinesZh-book-2021,Robinson-book-99}.
By the Smale Spectral Decomposition Theorem \cite{Smale67}, the non-wandering set of an Axiom A dynamical system can be decomposed into a disjoint union of closed, invariant, and transitive sets, which are known as
\textit{basic sets}.
Basic sets are classified into three types: attractors, repellers, and saddle-type.
A basic set is defined as \textit{trivial} if it consists of a single isolated trajectory; otherwise, it is considered \textit{nontrivial}.

This paper focuses on structurally stable, nonsingular, chaotic flows on closed orientable 3-manifolds whose basic sets consist solely of attractors and repellers.
For brevity, such flows will be referred to as being of the attractor-repeller type.
Bowen \cite{Bow70,Bow71} proved that an Axiom A dynamical system has positive topological entropy (and is thus chaotic) if and only if it possesses at least one nontrivial basic set.
Since any nontrivial one-dimensional basic set is of the saddle type, it follows that for a chaotic, nonsingular, attractor-repeller 3-flow, any nontrivial basic set must be two-dimensional.
Such a basic set is either an expanding attractor or a contracting repeller \cite{MedvedevZh2020-translation,MedvedevZh2022}.
The precise definitions for these terms are provided in Section \ref{s:prev}.

Let $\mathbb{T}^n$ denote the $n$-torus for $n \geq 2$.
By applying Smale's surgery operation to any Anosov automorphism on $\mathbb{T}^2$, one can construct what is known as a DA-diffeomorphism, $f: \mathbb{T}^2 \to \mathbb{T}^2$.
This is a structurally stable diffeomorphism whose non-wandering set consists of a one-dimensional expanding attractor and a source.
Consequently, the dynamical suspension $sus^t(f)$ over $f$ yields a structurally stable, chaotic, attractor-repeller type flow, which possesses a two-dimensional nontrivial attractor and a one-dimensional trivial repeller.
This construction demons\-t\-rates that the set of structurally stable, chaotic, nonsingular 3-flows of the attractor-repeller type is nonempty.
A natural question thus arises: which closed 3-manifolds admit such flows?

A foundational example of a chaotic, structurally stable, nonsingular 3-flow of the attractor-repeller type was constructed by Franks and Williams \cite{FranksWilliams1980}.
They constructed a nontransitive Anosov 3-flow whose non-wandering set consists of a two-dimensional expanding attractor and a two-dimensional contracting repeller.  Many examples were constructed in \cite{BeguinBonattiBinYu2017,JiagangYangBinYu-2022,Christy1993}.

Our first main result generalizes the work of \cite{MedvedevZh2022}, where it was proven that any closed, orientable 3-manifold $M$ admits an A-flow with an expanding attractor.
Here, we extend this result by providing an exact description of the non-wandering set for a structurally stable flow.
Throughout this paper, a \textbf{repelling trajectory} is defined as a hyperbolic periodic trajectory with a three-dimensional unstable manifold.
Such a trajectory constitutes a trivial repelling basic set.
We use the notation $NW(f^t)$ to denote the non-wandering set of a flow $f^t$.

\begin{theorem}\label{thm:any-three-man-admit-attractor-repeller}
Any closed, orientable 3-manifold $M$ admits a structurally stable, nonsingular flow $f^t$ whose non-wandering set, $NW(f^t)$, consists of a two-dimensional expanding attractor and a finite number of repelling periodic trajectories.
\end{theorem}

This result contrasts sharply with the case for diffeomorphisms. For instance, it has been proven that if a structurally stable diffeomorphism $f: M \to M$ on a closed 3-manifold $M$ (not necessarily orientable) possesses a two-dimensional expanding attractor, then $M$ must be the 3-torus $\mathbb{T}^3$ \cite{GinesZhuzhoma2005}.

In our second main result, we specify the set of repelling periodic trajectories from Theorem \ref{thm:any-three-man-admit-attractor-repeller} for the specific case where $M = \mathbb{S}^3$.

\begin{theorem}\label{thm:realization-for-link}
For any link $\{l_1,\ldots,l_k\} \subset \mathbb{S}^3$ that contains the figure-eight knot, there exists a structurally stable, nonsingular flow $f^t$ on $\mathbb{S}^3$ whose non-wandering set, $NW(f^t)$, consists of a two-dimensional expanding attractor and a finite number of repelling periodic trajectories whose union forms the given link $\{l_1,\ldots,l_k\}$.
\end{theorem}

Our third main result addresses the case of a flow on $\mathbb{S}^3$ whose non-wandering set consists of a single repelling periodic trajectory and a single expanding attractor.

\begin{theorem}\label{thm:conj-for-s-three-true}
Let $f^t$ be a structurally stable, nonsingular flow on $\mathbb{S}^3$ whose non-wandering set consists of a unique repelling periodic trajectory $l$ and an expanding attractor. Then $l$ cannot be a torus knot.
\end{theorem}

If we replace  $\mathbb{S}^3$ by general \(3\)-manifold, we have  the following weaker result:
\begin{theorem}\label{thm:revised-version}
Let $f^t$ be a structurally stable, nonsingular flow on a closed, orientable 3-manifold $M$ whose non-wandering set consists of a unique repelling periodic trajectory $l$ and an expanding attractor. Then $l$ is a nontrivial knot (i.e., there is no embedded disk $D^2\subset M$ with $\partial D^2=l$).
\end{theorem}

Several related results provide context for our work. Christy \cite{Christy1993} classified the restrictions of A-flows ($f^t|_{\Lambda_a}$) to their expanding attractors ($\Lambda_a$) up to an equivalence relation. However, it is important to note that a classification of these restrictions does not imply a classification of the basins of the attractors \cite{RobinsonWilliams76,ZhuzhIsaenkova2009}. In other work, B\'{e}guin, Bonatti, and Bin Yu \cite{BeguinBonattiBinYu2017} constructed nontransitive Anosov flows with two-dimensional attractors that have a prescribed entrance foliation, including some incoherent attractors. More recently, Bin Yu and Jiagang Yang \cite{JiagangYangBinYu-2022} classified expanding attractors on the figure-eight knot complement and nontransitive Anosov flows on the so-called Franks-Williams manifold. They proved that, up to topological equivalence, there exists a unique flow with an expanding attractor in the complement of the figure-eight knot on the 3-sphere and described the topological structure of this flow.

The remainder of this paper is organized as follows. Section \ref{s:prev} introduces the necessary definitions and reviews previous results. In Section \ref{s:proofs}, we prove our main theorems. Finally, Section \ref{s:conj-concl} presents several conjectures and discusses related open problems.

\textit{Acknowledgments}. This work was supported by the HSE University ``International Academic Cooperation" project. The research for this article began while Evgeny Zhuzhoma was visiting the School of Mathematical Sciences, Key Laboratory of Intelligent Computing and Applications, at Tongji University in Shanghai, China. He is grateful to Tongji University for its hospitality.

\section{Basic definitions and previous results}\label{s:prev}

\textsl{Structural stability}.
Let $F(M^n)$ denote the space of $C^1$-flows on a closed manifold $M^n$. This space is endowed with the strong $C^1$-topology. Two flows $f^t_1, f^t_2 \in F(M^n)$ are considered \textit{topologically equivalent} if there exists a homeomorphism $\varphi: M^n \to M^n$ that maps the trajectories of $f^t_1$ to the trajectories of $f^t_2$. A flow $f^t \in F(M^n)$ is \textit{structurally stable} if it has a neighborhood $V(f^t) \subset F(M^n)$ such that any flow $g^t \in V(f^t)$ is topologically equivalent to $f^t$.

\medskip
\noindent
\textsl{A-flows}.
According to Hayashi \cite{Hayashi97}, any structurally stable flow satisfies Smale's Axiom A. This means the non-wandering set is hyperbolic and the set of periodic trajectories is dense within it. Such flows are often referred to as A-flows.

For a smooth flow $f^t$ on a closed $n$-manifold $M^n$, where $n\geq 3$, an invariant, nonsingular set $\Lambda\subset M$ is defined as \textit{hyperbolic} if the tangent bundle over $\Lambda$, $T_{\Lambda}M$, admits a continuous, $Df^t$-invariant splitting into sub-bundles: $T_{\Lambda}M = E^{ss}_{\Lambda}\oplus E^t_{\Lambda}\oplus E^{uu}_{\Lambda}$. This splitting must satisfy the following conditions:
\begin{enumerate}
    \item The dimensions of the sub-bundles sum to $n$: $\dim E^{ss}_{\Lambda} + \dim E^t_{\Lambda} + \dim E^{uu}_{\Lambda} = n$.
    \item $E^t_{\Lambda}$ is the one-dimensional sub-bundle tangent to the flow trajectories.
    \item There exist constants $C_s > 0$, $C_u > 0$, and $\lambda \in (0,1)$ such that the stable and unstable sub-bundles exhibit exponential contraction and expansion:
    $$\Vert df^t(v)\Vert \leq C_s\lambda ^t\Vert v\Vert \text{ for } v\in E^{ss}_{\Lambda}, \quad \text{and} \quad \Vert df^{-t}(v)\Vert \leq C_u\lambda ^t\Vert v\Vert \text{ for } v\in E^{uu}_{\Lambda}, \quad \text{where } t>0.$$
\end{enumerate}

By the Smale Spectral Decomposition Theorem \cite{Smale67}, the non-wandering set $NW(f^t)$ of an A-flow $f^t$ decomposes into a disjoint union of closed, invariant, and transitive sets known as \textit{basic sets}. Furthermore, the hyperbolic structure of $NW(f^t)$ implies that for each point $x \in NW(f^t)$, there exist \textit{unstable}, \textit{stable}, \textit{strongly unstable}, and \textit{strongly stable manifolds}, denoted by $W^u(x)$, $W^s(x)$, $W^{uu}(x)$, and $W^{ss}(x)$, respectively.

Following Williams \cite{Williams74}, a basic set $\Omega$ is called an \textit{expanding attractor} if it is an attractor and its topological dimension equals the dimension of the unstable manifold $W^u(x)$ for any point $x \in \Omega$ (see the survey \cite{GrinesZh2006}). It was shown in \cite{MedvedevZh2022} (see also \cite{MedvedevZh2020-translation}) that an attractor $\Lambda_a$ of a 3-flow is expanding if and only if it is two-dimensional. In this case, $\Lambda_a$ is locally homeomorphic to the product of the Euclidean plane $\mathbb{R}^2$ and a Cantor set.

The following statement, proven in \cite{JiagangYangBinYu-2022} (see also \cite{BeguinBonattiBinYu2017}), is provided for reference. A similar result for Anosov flows was established earlier by Brunella \cite{Brunella1993}.
\begin{lemma}\label{lm:nghd-attractor-boundary-tori}
Let $f^t$ be an A-flow on a closed 3-manifold $M$, and let $\Lambda$ be an expanding attractor of $f^t$. Then there exists a connected attracting neighborhood $U$ of $\Lambda$ such that:
\begin{itemize}
    \item The boundary $\partial U$ is transverse to the flow.
    \item Each component of $\partial U$ is either a Klein bottle or a 2-torus. Furthermore, if $M$ is orientable, then every component of $\partial U$ is a 2-torus.
\end{itemize}
\end{lemma}

The following proposition is a crucial step in the proof of Theorem \ref{thm:any-three-man-admit-attractor-repeller}.

\begin{proposition}\label{prop:attractor-is-unique}
Let $f^t$ be an A-flow on a closed, orientable 3-manifold $M$. Suppose the non-wandering set of $f^t$ consists of a finite number of isolated repelling periodic trajectories, $l_1, \ldots, l_m$, and a set of expanding attractors, $\Lambda_1, \ldots, \Lambda_q$. Then there is only one such attractor; that is, $q=1$.
\end{proposition}

\textsl{Proof}. Each repelling periodic trajectory $l_i$ has a tubular neighborhood $U(l_i)$ that is homeomor\-p\-hic to a solid torus. Consequently, its boundary $\partial U(l_i)$ is a 2-torus. Since each $l_i$ is a hyperbolic repeller, we can choose these neighborhoods such that their boundaries are transverse to the flow $f^t$ for all $i=1,\ldots,m$.

Furthermore, each attractor $\Lambda_j$ is a connected set because it is the topological closure of a transitive trajectory. By Lemma \ref{lm:nghd-attractor-boundary-tori}, we can select disjoint, connected, attracting neighborhoods $U(\Lambda_1), \ldots, U(\Lambda_q)$ for the attractors $\Lambda_1, \ldots, \Lambda_q$ that satisfy the conditions of the lemma.

Since the non-wandering set consists only of these attractors and repellers, any trajectory escaping a neighborhood $U(l_i)$ must enter one of the attracting neighborhoods $U(\Lambda_j)$. This dynamic partitions the entire manifold $M$ into the disjoint basins of attraction for each $\Lambda_j$. As $M$ is a connected space, it cannot be decomposed into a disjoint union of multiple non-empty open sets (the basins). Therefore, there can be only one such basin, which implies there is only one attractor. Thus, $q=1$.
$\Box$

\medskip
\noindent
\textsl{Branched coverings}.
Let $M$ and $\widetilde{M}^3$ be closed, orientable 3-manifolds, and let $l \subset M$ be a knot (a smoothly embedded circle). For simplicity, we will provide the definition of a branched covering map over a knot; the definition over a link is analogous.

Let $\mathbb{D}^2 = \{z \in \mathbb{C} \mid |z| \leq 1\}$ be the unit disk in the complex plane $\mathbb{C}$. A \textit{dihedral group} $D_{2n}$ is a finite group of order $2n$ with the presentation $\{a,b \mid a^2=b^2=(ab)^n=1\}$. Geometrically, $D_{2n}$ can be represented as the group of symmetries of a regular $n$-polygon, $P_n$, inscribed in the unit circle $\partial\mathbb{D}^2$. This group consists of $n$ rotations and $n$ reflections.\footnote{The group $D_2$ is isomorphic to $\mathbb{Z}_2$, while $D_4$ is isomorphic to $\mathbb{Z}_2 \times \mathbb{Z}_2$.} For example, a rotation by $2\pi/n$ is given by the map $z \mapsto e^{i2\pi/n}z$, and a reflection across the real axis is given by the map $z \mapsto \bar{z}$.

We say that $p: \widetilde{M}^3\to M$ is an $m$-fold \textit{branched covering map} branched over a link $l$ if the restriction
$$ p|_{\widetilde{M}^3\setminus p^{-1}(l)}: \widetilde{M}^3\setminus p^{-1}(l)\to M\setminus l $$
is a standard (non-branched) $m$-fold covering map, where the preimage $p^{-1}(l)\subset\widetilde{M}^3$ is also a link with components $l_1, \ldots, l_k$ (pairwise disjoint embedded circles in $\widetilde{M}^3$) satisfying the following properties:
\begin{itemize}
    \item For a tubular neighborhood $U(l) \cong \mathbb{S}^1\times\mathbb{D}^2$ of $l$, its preimage $p^{-1}(U(l))$ is a disjoint union of tubular neighborhoods $U(l_1), \ldots, U(l_k)$ of $l_1, \ldots, l_k$, respectively, where each component $U(l_i)$ is also diffeomorphic to $\mathbb{S}^1\times\mathbb{D}^2$ for $i=1,\ldots,k$.
    \item For each $l_i$, there is a subgroup $G_i$ of the dihedral group such that the restriction
    $$p|_{U(l_i)}: U(l_i) \to U(l)$$
    is conjugate to an action of $G_i$.
    \item The preimage $p^{-1}(x)$ consists of $m$ points for every $x\in M\setminus l$.
\end{itemize}
If each subgroup $G_i$ consists only of rotation symmetries, then $p$ is called a \textit{cyclic} branched covering map.

The knot $l$ is called the \textit{branch set}. Similarly, one can define a branched covering map where the branch set is any link in $M$. For any given knot $l\subset M$ and any integer $m\geq 2$, there exists a manifold $\widetilde{M}^3$ and an $m$-fold branched covering map with branch set $l$ \cite{Rolfsen-book}.

A flow $\tilde{f}^t$ on $\widetilde{M}^3$ is called a \textit{lift} of the flow $f^t$ on $M$ (or a \textit{cover flow}) if the covering map $p$ takes every trajectory of $\tilde{f}^t$ onto a trajectory of $f^t$ while preserving the time direction. It is well-known that if $p: \widetilde{M}^3\to M$ is a standard (non-branched) covering map, any flow $f^t$ on $M$ can be lifted to a cover flow $\tilde{f}^t$ on $\widetilde{M}^3$; this is a direct consequence of the path-lifting property. For branched coverings, it follows directly from the definition that for a lift to exist, the branch set must be an invariant set under the flow $f^t$. The following proposition is necessary for the proof of Theorem \ref{thm:any-three-man-admit-attractor-repeller}.

\begin{proposition}\label{prop:cover-flow-for-branch-covering}
Let $f^t$ be a flow on $M$ with an isolated, repelling periodic trajectory $l$. Suppose $p: \widetilde{M}^3\to M$ is an $m$-fold branched covering map with branch set $l$. Then there exists a cover flow $\tilde{f}^t$ on $\widetilde{M}^3$ such that the preimage $p^{-1}(l)$ consists of isolated, repelling periodic trajectories for $\tilde{f}^t$.
\end{proposition}

\textsl{Proof}. We use the notation above. Note that the covering manifold $\widetilde{M}^3$ does not depend on the choice of the neighborhood $U(l)$ (in fact, it does not even depend on the Seifert surface for $l$). Since $l$ is a hyperbolic, repelling trajectory, by the Grobman-Hartman Theorem \cite{Robinson-book-99}, we can choose a neighborhood $U(l)\cong\mathbb{S}^1\times\mathbb{D}^2$ such that the restricted flow $f^t|_{U(l)}$ is invariant under any rotation of the form $(t,z)\to (t,z\cdot e^{i\phi})$. As a consequence, the flow $f^t|_{U(l)}$ is invariant under the central symmetry $(t,z)\to (t,-z)$ and any reflection. This symmetry implies the existence of a cover flow on each neighborhood $U(l_i)$ of $l_i$, where $p^{-1}(l)=\{l_1, \ldots, l_k\}$. On the complement of $\cup_{i=1}^kU(l_i)$, the map $p$ is a standard (non-branched) covering, so a lift also exists there. Patching these local lifts together yields the global cover flow $\tilde{f}^t$ on $\widetilde{M}^3$. By construction, each component $l_i$ is a periodic trajectory for $\tilde{f}^t$, and its repelling property is inherited from $l$. \
$\Box$

\section{Proofs of main results}\label{s:proofs}

The starting point for the proofs of Theorems \ref{thm:any-three-man-admit-attractor-repeller} and \ref{thm:realization-for-link} is the following \textit{canonical}, structurally stable, nonsingular flow $f_0^t$ on $\mathbb{S}^3$. Let $f: \mathbb{T}^2\to\mathbb{T}^2$ be a DA-diffeomorphism obtained from the Anosov diffeomorphism
$A=\begin{pmatrix}
    2 & 1 \\
    1 & 1
\end{pmatrix}: \mathbb{T}^2\to\mathbb{T}^2$.
Recall that $f$ is obtained by replacing the saddle fixed point $(0,0)$ of $A$ with three fixed points: two saddles and a source, denoted $O$. This procedure is known as Smale's surgery (see \cite{Robinson-book-99} for details).
As a result, $f$ is a structurally stable diffeomorphism whose non-wandering set consists of a one-dimensional expanding attractor, denoted $\lambda$, and the source $O$. The dynamical suspension of $f$, $sus(f)$, is then a nonsingular, structurally stable flow on the mapping torus $M_f=\mathbb{T}^2\times [0,1]/(x,1)\sim (f(x),0)$.
The non-wandering set of $sus(f)$ consists of an isolated, repelling periodic trajectory $l_0$ (generated by $O$) and a two-dimensional expanding attractor $\Lambda_a$ (generated by $\lambda$).

Let $l$ be a knot in a 3-manifold $N^3$, and let $T(l)$ be its tubular neighborhood. We remove $T(l)$ from $N^3$ and then glue it back using a diffeomorphism
$\phi: \partial T(l)\to\partial\left(N^3\setminus T(l)\right)$. This yields a new manifold $T(l)\cup_{\phi}\left(N^3\setminus T(l)\right)=N^3_l$. This reconstruction $N^3\to N^3_l$ is called a \textit{Dehn surgery} along $l$.

According to \cite{Christy1993,Thurston1997-book} (see also \cite{Francis1983}), there exists a Dehn surgery along $l_0$ such that this surgery transforms the mapping torus
$M_f=\mathbb{T}^2\times [0,1]/(x,1)\sim (f(x),0)$ into the 3-sphere $\mathbb{S}^3$ and transforms the flow $sus(f)$ into a flow $f_0^t$ on $\mathbb{S}^3$, whose non-wandering set $NW(f_0^t)$ consists of an isolated, repelling periodic trajectory and an expanding attractor. We will henceforth denote the repelling periodic trajectory by $l_1$ and the expanding attractor by $\Lambda_a$. It is well-known that $l_1$ is the figure-eight knot \cite{JiagangYangBinYu-2022}.

\medskip
\textsl{Proof of Theorem \ref{thm:any-three-man-admit-attractor-repeller}}. Consider the canonical structurally stable flow $f_0^t$ on $\mathbb{S}^3$ whose non-wandering set is $NW(f_0^t)=l_1\cup\Lambda_a$, where $\Lambda_a$ is an expanding attractor and $l_1$ is a repelling periodic trajectory that is the figure-eight knot. Since the figure-eight knot $l_1$ is a universal knot \cite{HildenLosanoMontesinos1985}, every closed, orientable 3-manifold $M$ can be represented as a covering of $\mathbb{S}^3$ branched over the link $l_1$, via a branched covering map $p: M\to\mathbb{S}^3$. According to Proposition \ref{prop:cover-flow-for-branch-covering}, there exists a flow $f^t_1$ on $M$ that is a lift of $f^t_0$. By construction, $f^t_1$ is a nonsingular flow.

Since $NW(f_0^t)=l_1\cup\Lambda_a$, the non-wandering set of $f^t_1$ is the preimage $p^{-1}(NW(f_0^t))$, which consists of the expanding attractor $p^{-1}(\Lambda_a)$ and the link $p^{-1}(l_1)$ (a set of finitely many isolated, repelling periodic trajectories). Thus, $f^t_1$ is a chaotic, nonsingular flow of the attractor-repeller type. By Proposition \ref{prop:attractor-is-unique}, $f^t_1$ has a unique expanding attractor. Moreover, the 3-dimensional unstable manifold of each isolated repelling periodic trajectory intersects the stable manifold of any trajectory from the attractor $p^{-1}(\Lambda_a)$ transversally. It follows from \cite{Hayashi97} that $f^t_1$ is structurally stable.
$\Box$

\medskip
\textsl{Proof of Theorem \ref{thm:realization-for-link}}. Set the given link $\mathcal{L}=\{l_1,\ldots,l_k\}$. Consider the canonical structurally stable flow $f_0^t$ on $\mathbb{S}^3$, whose non-wandering set is $NW(f_0^t)=l_1\cup\Lambda_a$, where $\Lambda_a$ is an expanding attractor and $l_1$ is a repelling periodic trajectory. If $\mathcal{L}$ consists of the single knot $l_1$, then there is nothing to prove, and we set $f^t=f_0^t$. Let us now suppose that $\mathcal{L}$ contains knots other than $l_1$, i.e., $k\geq 2$.

According to \cite{Ghrist1997} (see also \cite{GhristHolmesSullivan}), there exist periodic trajectories $l_2, \ldots, l_k$ of $f_0^t$ such that the set $\{l_1, l_2, \ldots, l_k\}$ forms the desired link $\mathcal{L}$. By applying Smale's surgery to each trajectory $l_2, \ldots, l_k$, we obtain a new flow $f^t$ in which $l_2, \ldots, l_k$ have become isolated, repelling periodic trajectories. Since Smale's surgery is a local bifurcation, the resulting flow $f^t$ on $\mathbb{S}^3$ has a non-wandering set, $NW(f^t)$, that consists of a (two-dimensional) expanding attractor and the finitely many repelling periodic trajectories that form the link $\mathcal{L}=\{l_1, \ldots, l_k\}$.
$\Box$

\medskip
\textsl{Proof of Theorem \ref{thm:conj-for-s-three-true}}. Suppose, for contradiction, that $l$ is a $(p,q)$-torus knot, and let $U(l)$ be a tubular neighborhood of $l$. Then $\mathbb S^3\setminus U(l)$ is the Seifert manifold
\[
N=(D^2;(p,1),(q,1)).
\]

Let $X$ denote the vector field induced by the flow $f^t$. The manifold $N$ contains an expanding attractor $\Lambda_a$ for $X$. Without loss of generality, we may assume that $X$ points outward along $\partial N$. Moreover, the stable foliations of $\Lambda_a$ can be taken transverse to $\partial N$, thereby inducing a 1-foliation $f^s$ on $\partial N$ (see \cite{BeguinBonattiBinYu2017,FranksWilliams1980}).

Now take a copy $N_1$ of $N$ equipped with the reversed vector field $-X$. Clearly, $-X$ has a contracting repeller $\Lambda_r$ homeomorphic to $\Lambda_a$. Its unstable foliation is transverse to $\partial N_1$, inducing a foliation $f^u$ on $\partial N_1$.

According to the definition of MS-foliations given in the preprint \cite{BeguinBonattiBinYu2025}
(see Summary of the preprint \cite{BeguinBonattiBinYu2025}), both $f^s$ on $\partial N$ and $f^u$ on $\partial N_1$ are MS-foliations.
Therefore, by Proposition \ref{appendix-A}, there exists a gluing map
\[
\phi:\partial N_1\to\partial N,
\]
isotopic to the identity, such that $f^s\pitchfork \phi(f^u)$. Applying Theorem 1.5 of \cite{BeguinBonattiBinYu2017}, we glue $N$ and $N_1$ along their boundaries under the isotopy class that we chose  to obtain a closed $3$–manifold $L$ supporting a non-transitive Anosov flow $\varphi_t$. Note that \(L\) is homeomorphic to  the manifold obtained by doubling $N$ along the boundaries, therefore \(L\) is homeomorphic to the Seifert manifold
\[
L\cong (\mathbb S^2;(p,1),(p,1),(q,1),(q,1)).
\]

According to \cite{Ghys-1984}, up to finite cover, the Anosov flow $\varphi_t$ is topologically equivalent to a geodesic flow; in particular, it is transitive. This contradicts to the fact that
\(\varphi_t\) is non-transitive. The proof of Theorem \ref{thm:conj-for-s-three-true}  is complete.
$\Box$

\medskip
\textsl{Proof of Theorem \ref{thm:revised-version}}. Suppose, for contradiction, that $l$ is a trivial knot, and let $U(l)$ be a tubular neighborhood of $l$. Denote
\[N \;=\; M \setminus U(l) \;\cong\; M \,\#\, (D^2 \times \mathbb{S}^1).\] Following the same argument as in the proof of Theorem~\ref{thm:conj-for-s-three-true}, take another copy $N_1$ of $N$ and glue $N_1$ to $N$ to obtain a closed $3$-manifold.\[W \;\cong\; M \,\#\, (\mathbb{S}^2 \times \mathbb{S}^1) \,\#\, M,\]which supports an Anosov flow. On one hand, the weak stable manifolds of the Anosov flow provides a \(2\)-foliation \(\mathcal{F}^s\). Notice that each leaf of \(\mathcal{F}^s\)
is homeomorphic to one of \(\mathbb{R}^2\), an open annulus, an open Mobius band, consequently, \(\mathcal{F}^s\)
does not contain any compact leaf, therefore, due to Lemma 4.24 of  \cite{Calegari-2007},
\(\mathcal{F}^s\) is a taut foliation.
 On the other hand, since $W$ is reducible, essentially due to Novikov theorem (see Theorem 4.35 of \cite{Calegari-2007}),  it cannot carry any taut foliation. This contradiction shows that $l$ cannot be trivial.
$\Box$

\section{Further discussions}\label{s:conj-concl}

Here we discuss some problems and formulate several conjectures concerning chaotic, nonsingular 3-flows of the attractor-repeller type.

According to \cite{HildenLosanoMontesinos1985}, there exists a 6-fold branched covering map $\mathbb{S}^3\to\mathbb{S}^3$ branched over the figure-eight knot $l$, whose preimage contains the Borromean rings. This result, combined with Proposition \ref{prop:cover-flow-for-branch-covering}, implies the existence of a structurally stable, nonsingular flow $f^t$ on $\mathbb{S}^3$ whose non-wandering set $NW(f^t)$ consists of an expanding attractor and finitely many repelling periodic trajectories. In light of Theorem \ref{thm:conj-for-s-three-true} and Theorem \ref{thm:revised-version}, it is natural to formulate the following conjecture.
\begin{conjecture}
Let $f^t$ be a structurally stable, nonsingular flow on $\mathbb{S}^3$ whose non-wandering set consists of an expanding attractor and a finite set of repelling periodic trajectories $\{l_1, \ldots, l_k\}$. Then this set of trajectories forms a nontrivial link, or at least one trajectory $l_j$ is a nontrivial knot.
\end{conjecture}

Note that this conjecture is not true for $\mathbb{S}^2\times\mathbb{S}^1$. Indeed, consider a structurally stable diffeomorphism $f: \mathbb{S}^2\to\mathbb{S}^2$ whose non-wandering set consists of a Plykin attractor and four isolated sources. The dynamical suspension of this diffeomorphism is a structurally stable, nonsin\-gular flow on $\mathbb{S}^2\times\mathbb{S}^1$ with an expanding attractor and four repelling periodic trajectories, all of which are unlinked, trivial knots.

Jiagang Yang and Bin Yu  proved a notable result regarding expanding attractors \cite{JiagangYangBinYu-2022}. They showed that, up to topological equivalence, there is a unique flow with an expanding attractor supported on the complement of the figure-eight knot. They denoted this flow by $Y^0_t$. It would be interesting to construct a flow topologically distinct from $Y^0_t$.

\begin{conjecture}
There exists a nonsingular flow on $\mathbb{S}^3$ whose non-wandering set consists of an expanding attractor and a unique repelling periodic trajectory that is not the figure-eight knot.
\end{conjecture}

Let $p: \widetilde{M}^3\to\mathbb{S}^3$ be a 2-fold cyclic branched covering map branched over the figure-eight knot $l$. Then the preimage $p^{-1}(l)$ is a knot in $\widetilde{M}^3$.
It is known that there is a nonsingular flow on $\widetilde{M}^3$, say $f^t_*$, whose non-wandering set consists of an expanding attractor and the repelling periodic trajectory $p^{-1}(l)$. In light of \cite{JiagangYangBinYu-2022}, it is interesting to consider the following conjecture.

\begin{conjecture}
Up to topological equivalence, there exists a unique flow on $\widetilde{M}^3$ (namely, the class of $f^t_*$) whose non-wandering set consists of an expanding attractor and a unique repelling periodic trajectory.
\end{conjecture} 

\clearpage

\appendix
\section{Summary of the preprint \cite{BeguinBonattiBinYu2025}}
In this appendix we briefly summarize certain definitions and results from
B\'{e}guin, Bonatti, and Yu \cite{BeguinBonattiBinYu2025}, since the paper is not yet published.
In particular, we recall the notion of an \emph{MS-foliation} on the torus $\mathbb{T}^2$
and state Proposition~2 of their work, which is used in the proof of the Theorem~\ref{thm:conj-for-s-three-true}.

\medskip
\noindent
\textsl{MS-foliation}. A foliation on a torus  $\mathbb{T}^2$ is called a \textsl{MS-foliation} if it satisfies the following properties:
\begin{enumerate}
    \item It has only finitely many compact leaves;
    \item Every half leaf approaches a compact leaf;
    \item Every compact leaf can be oriented such that its holonomy is a contraction.
\end{enumerate}

\begin{proposition}\label{appendix-A}
(Proposition~2 in \cite{BeguinBonattiBinYu2025})
Let $f$ be an MS-foliation on $\mathbb{T}^2$.
Then there exists a self-homeomorphism $\varphi$ of $\mathbb{T}^2$, isotopic to the identity, such that $f \pitchfork \varphi(f)$.
\end{proposition}

\textsl{Proof}.
Let $c^0_1,\dots,c^0_m$ be the compact leaves of $f$ such that there is no compact leaf of $f$ between $c^0_i$ and $c^0_{i+1}$.
Fix an orientation on every compact leaf of $f$ so that $c^0_i$ and $c^0_j$ are co-oriented.
We call this orientation the \emph{topological orientation} of $c^0_i$.
The \emph{dynamical orientation} of $c^0_i$ is defined by the contracting holonomy direction of $f$ at $c^0_i$.
If the dynamical orientation and the topological orientation of $c^0_i$ coincide, we call $c^0_i$ positive; otherwise, we call $c^0_i$ negative.

Let $\delta_0=(c^0_1,\dots,c^0_m)$ be the compact-leaf cyclic sequence of $f$.
Take a copy of $f$, denoted $g$, with compact-leaf cyclic sequence $\delta_1=(c^1_1,\dots,c^1_m)$ such that $c^1_i$ corresponds to $c^0_i$.
Next, isotope $g$ according to the following rules:
\begin{itemize}
\item If $c^0_i$ is positive, then $c^1_i$ is placed between $c^0_i$ and $c^0_{i+1}$;
\item If $c^0_i$ is negative, then $c^1_i$ is placed between $c^0_{i-1}$ and $c^0_i$;
\item $c^1_i$ and $c^0_i$ are adjacent.
\end{itemize}

It is straightforward to check that $f \pitchfork g$, and that $f$ and $g$ are leaf-conjugate via a map $\varphi$ isotopic to the identity.
$\Box$

\bigskip
\noindent
School of Mathematical Sciences, Key Laboratory of Intelligent Computing and Applications(Ministry of Education), Tongji University, Shanghai, China

\noindent
\textit{E-mail:} binyu1980@gmail.com,  zhentao\_lai@126.com

\bigskip
\noindent
National Research University Higher School of Economics, Nizhny Novgorod, Russia

\noindent
\textit{E-mails:} medvedev-1942@mail.ru, zhuzhoma@mail.ru

\end{document}